\documentclass{article}
\usepackage[colorlinks=true, linkcolor=blue, citecolor=blue,urlcolor=black]{hyperref}
\usepackage[utf8]{inputenc} % allow utf-8 input
\usepackage[T1]{fontenc}    % use 8-bit T1 fonts
\usepackage{hyperref}       % hyperlinks
\usepackage{url}            % simple URL typesetting
\usepackage{booktabs}       % professional-quality tables
\usepackage{amsfonts}       % blackboard math symbols
\usepackage{nicefrac}       % compact symbols for 1/2, etc.
\usepackage{microtype}      % microtypography
\usepackage{xcolor}         % colors
\usepackage{graphicx}
\usepackage{forest}
\usepackage{algorithm}
\usepackage[noend]{algpseudocode}
\usepackage[utf8]{inputenc}
\usepackage{comment}
\usepackage{amsmath}
\usepackage{amssymb}
\usepackage{amsthm}
\usepackage{geometry}
\usepackage{graphicx}
\graphicspath{ {./images/} }
\usepackage{amsmath,amssymb}
\usepackage{fullpage}
\usepackage{color,soul}   %May be necessary if you want to color links and highlight stuff
\usepackage{hyperref}
\usepackage{todonotes}
\hypersetup{
	colorlinks=true, %set true if you want colored links
	linktoc=all,     %set to all if you want both sections and subsections linked
	linkcolor=blue,  %choose some color if you want links to stand out
}

\newtheorem{theorem}{Theorem}
\newtheorem{example}{Example}
\newtheorem{definition}{Definition}
\newtheorem{conjecture}{Conjecture}

\title{Stochastic Couplings and Bijections from the Symmetric Group to Itself}

\author{
	William Chang \\
	University of Southern California\\
	\texttt{chan087@usc.edu}
}
\begin{document}
	
	\maketitle
	
\begin{abstract}
Inspired by the Stochastic processes described by the Feller Coupling and Chinese Restaurant Processes, we create four different bijections from words in the set $[1]\times [2] \times\cdot \times[n]$ to $S_n$. We then compose these maps with their inverse to obtain a total of six bijections $S_n \to S_n$. Following that, we investigate the fixed points ($1$-cycle) and higher $k$-cycles of these maps. We characterized some of their properties completely as well as empirically showing the complexity of the higher $k$-cycle structures for these maps. 
\end{abstract}
\section{Introduction}

\indent \indent The Feller Coupling and Chinese Restaurant processes are both stochastic processes which generate permutations $S_n$ of $[n] := \{1,2,...,n\}$. The Feller Coupling was motivated by Feller \cite{feller1945fundamental}, while the Chinese Restaurant was motivated by Aldous \cite{aldous1985exchangeability}. Both of these were first explicitly described in Arratia \cite{arratia1992poisson}. The Chinese Restaurant process is as follows. Define a sequence of random variables $X_1, X_2,...,X_n$ with probabilities,

\begin{equation}
    P(X_i = j) = \begin{cases}
    \frac{\theta}{\theta+i-1}, \quad j = i\\
    \frac{1}{\theta + i - 1}, \quad j = 1,2,...,i-1
    \end{cases}
\end{equation}

We can now use this sequence of random variable to generate the cycles of a permutation. We start the cycle with $1$. For every integer $i\geq 2$, $i$ is placed to the immediate right of $j\leq i -1$ in the same cycle with probability $P(X_i = j)$, or it starts a new cycle with probability $P(X_i = i)$. On the other hand, the Feller Coupling is considers a sequence of independent Bernoulli random variables $Y_1, Y_2,... Y_n$ with probabilities given by 
\begin{equation}
    \mathbb{P}(Y_i = 1) = \frac{\theta}{\theta + j - 1}, \quad j \in [n]
\end{equation}

We can now construct a permutation in cycle notation by first starting with $1$, and for $i \geq 2$, if $Y_i= 1$ we start a new cycle putting the smallest unused integer in that cycle. If $Y_i= 0$, then we pick one of the remaining $n - i+1$ integers uniformly at random and place it to the immediate right of $1$ in the same cycle.  

There have been many works on both of these couplings. On one hand, Feller Coupling is often studied with the Ewen's Distribution with parameter $\theta$ introduced in Ewens \cite{ewens1972sampling} defined as follows
\begin{equation}
    P(a_1,...,a_n;\theta) = \frac{n!}{\theta(\theta+1)\cdots (\theta + n - 1)}\prod_{j=1}^n \frac{\theta^{a_j}}{j^{a_j}a_j!}
\end{equation}

where $\sum_{i=1}^k ia_i = n$. This formula describes the probability that in a random sample of $n$ gametes classified according to the gene at a specific locus there are  $a_i$ alleles represented $i$ times in the sample. $\theta$ represents the population mutation rate. In \cite{da2021feller} they study the Feller Coupling for random derangement of $[n]$ under the Ewens distribution with parameter $\theta$ arising as the weak limit as $n \to \infty$. Then, in \cite{ignatov1981point}, they describe a nice construction of said permutation in the case that $\theta =1$, and finally, \cite{najnudel2020feller} extends this approach for general $\theta > 0$. 

The Chinese Restaurant process has been used in various works related to modeling topic hierarchies (the study of taxonomies) such as Blei \cite{blei2010nested}, Griffiths \cite{griffiths2003hierarchical}, Wang \cite{wang2009variational}, and Kim \cite{kim2012modeling}. These have been used in areas such as taxonomies of images in works such as Bart \cite{bart2008unsupervised} and Sivic \cite{sivic2008unsupervised}. There are also applications to different mixture models such as Gaussian mixture models West\cite{west1993hierarchical}, hidden Markov models Beal \cite{beal2001infinite}, and mixtures of experts Rasmusen \cite{rasmussen2001infinite}.  

In section \ref{section:maps} we consider words belonging to $[1] \times [2] \times \cdots \times [n]$. It's clear that there are $n!$ such possible words, which allows us to create bijections from these words to permutations of $[n]$. We develop a total of $4$ types of bijections $f_{1,n}$, $f_{2,n}$, $f_{3,n}$, and $f_{4,n}$ from $[1] \times [2] \times \cdots \times [n]$ to $S_n$.  $f_{1,n}$ and $f_{2,n}$ maps words to permutations in one line notation, while $f_{3,n}$ and $f_{4,n}$ maps words to permutations in cycle notation. $f_{1,n}$ and $f_{3,n}$ are motivated by the Chinese Restaurant Process, while $f_{2,n}$ and $f_{4, n}$ are motivated by the Feller Coupling. In section \ref{section:bijection}, we compose these maps and their inverse to obtain bijections $S_n \to S_n$. Namely, for $i,j \in \{1, 2, 3, 4\}$ we compose $f_{i, n}$ with $f_{j, n}^{-1}$ to get bijections $S_n \to S_n$. The goal of this paper is to investigate the fixed points ($1$-cycle) and higher $k$-cycles of these maps. We completely characterize $f_{1,n}\circ f_{2,n}^{-1}$ by showing that this maps permutations to their inverse permutation but written backwards. Building on this map we provide a small modification that maps permutations to their inverses, and also develop a new sequence in Appendix \ref{appendix:f12}. For the maps $f_{1,n}\circ f_{3,n}^{-1}$ and $f_{2,n}\circ f_{4,n}^{-1}$ we characterize the fixed points (1-cycles) completely. For all the maps we also run numerical simulations in Appendix \ref{appendx:numerical}, empirically showing that there exists high order cycle structures even for small values of $n$.  

\section{Bijections $[1]\times [2]\times \cdots \times [n] \to S_n$}\label{section:maps}
To make describing the maps $f_{1,n}$, $f_{2,n}$, $f_{3,n}$, and $f_{4,n}$ easier, let us start with a definition.
\begin{definition}
Let $w \in [1]\times [2] \times \cdots \times [n]$ be a word. We denote the $i$th position counting from the left of a word $w$ or permutation $\sigma$ as $w[i]$ and $\sigma[i]$ respectively. Thus, if $w = 121$ then $w[1] = 1, w[2] = 2, w[3] = 1$. Another example is if $\sigma = 213\in S_3$, then $\sigma[1] = 2, \sigma[2] = 1, \sigma[3] = 3$.
\end{definition}

Let us define $f_{1,n}: [1]\times \cdots \times [n] \to S_n$ as follows.
\begin{definition}
 $f_{1,n}(w):[1] \times \cdots \times [n] \to S_n$ is $f_{1,n}(w')$ with $n$ inserted at the $w[n]$th position from the right. Here, $w'$ is $w$ but with only the first $n-1$ characters. The base case is $f_{1,n}(1) = 1$. 
\end{definition}

\begin{example}
For example suppose $n = 3$ and $w = 113$. Since $w[2] = 1$, we place $2$ in the rightmost spot giving us the permutation $12$ then since $w[3] = 3$ we place $3$ in position $3$ in the $3$rd spot counting from the right which leave us with the permutation $312$. Figure \ref{tree_f_1,n} has the pictorial representation of $f_{1,n}$ for $n = 3$. 

\end{example}

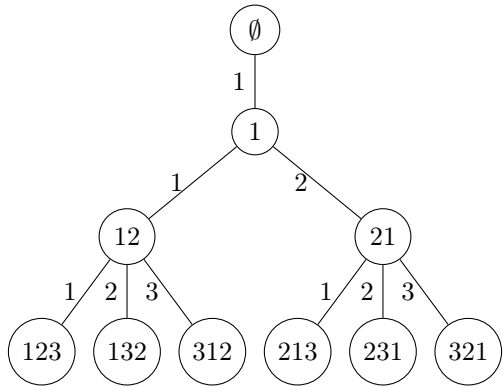
\begin{figure}
\centering
\begin{forest}
for tree={circle,draw, l sep=20pt}
[$\emptyset$ 
   [1, edge label={node[midway,left] {1}}
   [12, edge label={node[midway,left] {1}}
      [123, edge label={node[midway,left] {1}} ] 
      [132, edge label={node[midway,left] {2}}] 
      [312, edge label={node[midway,left] {3}}]
    ]
    [21, edge label={node[midway,left] {2}}
      [213, edge label={node[midway,left] {1}}] 
      [231, edge label={node[midway,left] {2}}] 
      [321, edge label={node[midway,left] {3}}]
    ]
  ]
]
\end{forest}
\caption{Pictorial representation of $f_{1,n}$. The edges are labeled with the numbers of our word $w$ with the top row of edges representing possible values for $w[1]$, second row of edges representing possible values for $w[2]$, etc. The circles show the process of building the permutation using the numbers on the edges.}
\label{tree_f_1,n}
\end{figure}

Let us now define our other map, $f_{2,n}$ 

\begin{definition}
$f_{2,n}: [1]\times \cdots \times [n] \to S_n$ is defined as follows. Let $w\in [1]\times \cdots \times [n]$ be expressed as $w = a_1...a_n$. This time we read the word backwards, starting from $w[n]$ and for each $i$ place in the rightmost spot the $w[i]$th smallest number that has yet to be in the permutation.
\end{definition}

We have the following example. 

\begin{example}
Suppose $n = 3$ and $w = 113$. Since $a_3 = 3$, we place the $3$rd smallest element in $\{1, 2, 3\}$ in the rightmost spot giving us $3$. Then since $w[2] = 1$ we place the first smallest element from $\{2, 3\}$ in the rightmost spot which leave us with the permutation $31$. Finally, since $a_1 = 1$ as always, we place the only remaining element in the very left leaving us with the permutation $312$. Figure \ref{tree_f_2,n} has the pictorial representation of $f_{2,n}$ when $n = 3$. 
\end{example}
\begin{figure}
 \centering  
\begin{forest}
for tree={circle,draw, l sep=20pt}
[$\emptyset$ 
   [1, edge label={node[midway,left] {1}}
   [12, edge label={node[midway,left] {1}}
      [123, edge label={node[midway,left] {1}} ] 
    ]
    [13, edge label={node[midway,left] {2}}
      [132, edge label={node[midway,left] {1}}] 
    ]
  ]
  [2, edge label={node[midway,left] {2}}
   [21, edge label={node[midway,left] {1}}
      [213, edge label={node[midway,left] {1}} ] 
    ]
    [23, edge label={node[midway,left] {2}}
      [231, edge label={node[midway,left] {1}}] 
    ]
  ]
  [3, edge label={node[midway,left] {3}}
   [31, edge label={node[midway,left] {1}}
      [312, edge label={node[midway,left] {1}} ] 
    ]
    [32, edge label={node[midway,left] {2}}
      [321, edge label={node[midway,left] {1}}] 
    ]
  ]
]
\end{forest}
 \caption{Tree diagram for $f_{2,n}$}
 \label{tree_f_2,n}
\end{figure}
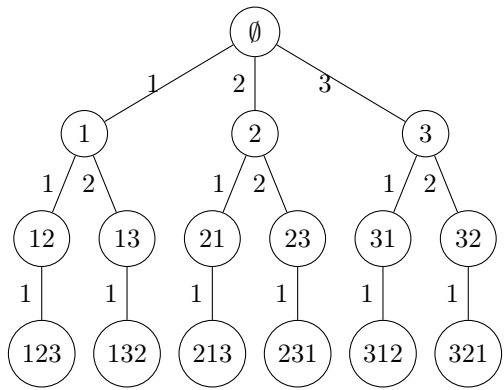

Let us now define $f_{3,n}: [1] \times \cdots \times [n] \to S_n$ which defines map from words to permutations in their cycle notation.

\begin{definition}
$f_{3,n}: [1] \times \cdots \times [n] \to S_n$ is defined as follows. For any $w \in [1] \times \cdots \times [n]$, we can build $f(w)$ by starting with $1$ and reading the word $w$ backwards. At each letter $w[i]$ of $w$, we insert a number in the cycle notation from the right of $f(w)$ as follows, 
\begin{itemize}
    \item If $w[i] > 1$, insert the $(w[i]-1)$-th smallest unused letter of $f(w)$ at this state. 
    
    \item If $w[i] =1$, close off the current cycle, and start a new one. For the new cycle put smallest unused letter of $f(w)$ as the first element of this cycle.  
\end{itemize}
\end{definition}

We have the following example to illustrate this map. 
\begin{example}
Suppose for $n=3$ we have $w = 113$. To evaluate $f(w)$ we start with the partial cycle $(1$. Then since $w[3] = 3>1$ we insert the 3rd smallest unused letter from the right. $1$ is always unused and as no other letters currently exist in $f(w)$, we insert $3$ from the right to get $(13$. Then since $w[2] = 2>1$ we look for the 2nd smallest unused letter. The only letter that is unused is $3$ so the second smallest unused letter is $2$. Inserting that after to the right gives $(132$. Finally $w[1] = 1$ so we start a new cycle put the smallest unused letter.  Figure \ref{tree_f_3,n} has the pictorial representation of $f_{3,n}$ when $n = 3$. 
\end{example}

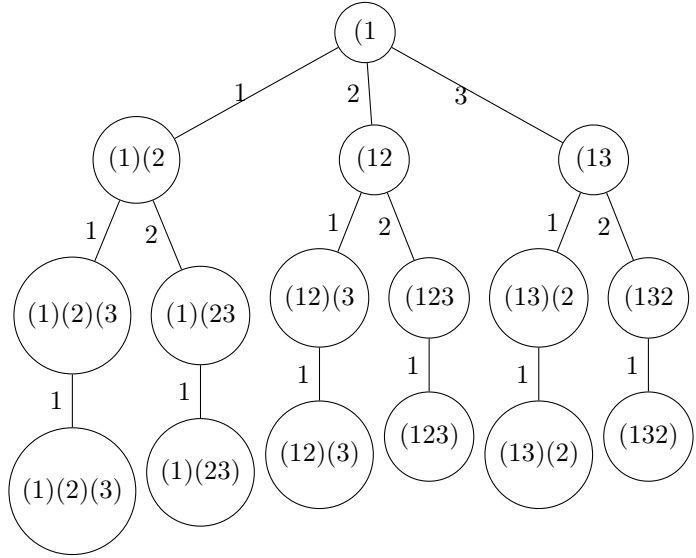
\begin{figure}\label{tree_f_3,n}
 \centering  
\begin{forest}
for tree={circle,draw, l sep=20pt}
[$(1$ 
   [(1)(2, edge label={node[midway,left] {1}}
   [(1)(2)(3, edge label={node[midway,left] {1}}
      [(1)(2)(3), edge label={node[midway,left] {1}} ] 
    ]
    [(1)(23, edge label={node[midway,left] {2}}
      [(1)(23), edge label={node[midway,left] {1}}] 
    ]
  ]
  [(12, edge label={node[midway,left] {2}}
   [(12)(3, edge label={node[midway,left] {1}}
      [(12)(3), edge label={node[midway,left] {1}} ] 
    ]
    [(123, edge label={node[midway,left] {2}}
      [(123), edge label={node[midway,left] {1}}] 
    ]
  ]
  [(13, edge label={node[midway,left] {3}}
   [(13)(2, edge label={node[midway,left] {1}}
      [(13)(2), edge label={node[midway,left] {1}} ] 
    ]
    [(132, edge label={node[midway,left] {2}}
      [(132), edge label={node[midway,left] {1}}] 
    ]
  ]
]
\end{forest}
 \caption{Tree diagram for $f_{3,n}$}
\end{figure}

Finally, we now define $f_{4,n}: [1] \times \cdots \times [n] \to S_n$ which also defines map from words to permutations in their cycle notation. 

\begin{definition}
$f_{4,n}: [1] \times \cdots \times [n] \to S_n$ is defined as follows. For any $w \in [1] \times \cdots \times [n]$ we can recursively build the permutation $f(w)$ in cycle notation with the base case $f_{3,n}(1) = (1)$. Let $w'$ be the first $n-1$ letters of $w'$. 
\begin{itemize}
    \item If $w[n] > 1$, insert $n$ after $w[n] - 1$ in $f(w')$.
    
    \item If $w[n] = 1$, start a new cycle in $f(w')$ with $n$ as the first letter.
\end{itemize}

\end{definition}

We have the following example to illustrate this map. 

\begin{example}
Suppose for $n=3$ we have $w = 113$. To evaluate $f(w)$ we start with the partial cycle $(1$. Then since $w[3] = 3>1$ we insert the 3rd smallest unused letter from the right. $1$ is always unused and as no other letters currently exist in $f(w)$, we insert $3$ from the right to get $(13$. Then since $w[2] = 2>1$ we look for the 2nd smallest unused letter. The only letter that is unused is $3$ so the second smallest unused letter is $2$. Inserting that after to the right gives $(132$. Finally $w[1] = 1$ so we start a new cycle put the smallest unused letter.  Figure \ref{tree_f_4,n} has the pictorial representation of $f_{4,n}$ when $n = 3$. 
\end{example}

\begin{figure}\label{tree_f_4,n}
\centering
\begin{forest}
for tree={circle,draw, l sep=20pt}
[$\emptyset$ 
   [(1), edge label={node[midway,left] {1}}
   [(1)(2), edge label={node[midway,left] {1}}
      [(1)(2)(3), edge label={node[midway,left] {1}} ] 
      [(13)(2), edge label={node[midway,left] {2}}] 
      [(1)(23), edge label={node[midway,left] {3}}]
    ]
    [(12), edge label={node[midway,left] {2}}
      [(12)(3), edge label={node[midway,left] {1}}] 
      [(132), edge label={node[midway,left] {2}}] 
      [(123), edge label={node[midway,left] {3}}]
    ]
  ]
]
\end{forest}
\caption{Tree Diagram for $f_{4,n}$}
\end{figure}
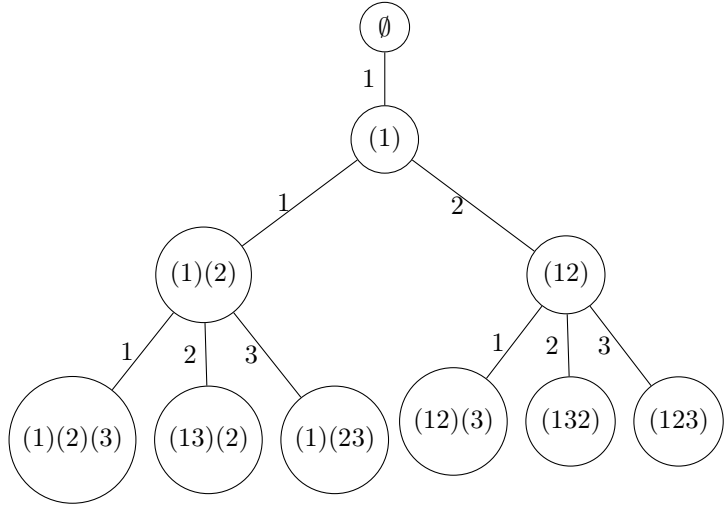

\section{Bijections $S_n \to S_n$}\label{section:bijection}
In this section, we will take $f_{1,n}$, $f_{2,n}$, $f_{3,n}$, and $f_{4,n}$ which maps from words $[1] \times [2] \times \cdots \times [n]$ to permutations and compose it with their inverses. This gives bijections $S_n \to S_n$. Namely, for $i<j \in \{1, 2, 3, 4\}$ we consider the six possible compositions $f_{i, n} \circ f_{j, n}^{-1}$. In the next subsection we will demonstrate this process with the map $f_{1,n}\circ f_{2,n}^{-1}$. We will then investigate the fixed points ($1$-cycle) and higher $k$-cycles of these maps. 

\subsection{The Bijection $f_{1,n}\circ f_{2,n}^{-1}$}

 For any value of $n$, $f_{1,n}$ and $f_{2,n}$ are both bijections between $[1]\times \cdots \times [n]$ and $S_n$, so we can consider $f_{1,n}\circ f_{2,n}^{-1}: S_n \to S_n$. Considering the trees from the Figures \ref{tree_f_1,n} and \ref{tree_f_2,n}, we can make list the words to their outputs under $f_{1,n}$ and $f_{2,n}$ like in
 Table \ref{table_f}, then the map $f_{1,n}\circ f_{2,n}^{-1}$ is simply the middle column of this table to the rightmost column. We are interested in the fixed points of this map and will prove in the next theorem that the number of fixed points is equal to the number of involution of $S_n$.

\begin{table}[]
        \centering
        \begin{tabular}{c|c|c}
            word & $f_{2,n}$ & $f_{1,n}$ \\
            \hline
             111 & 123 & 123 \\
             112 & 213 & 132 \\
             113 & 312 & 312 \\
             121 & 132 & 213 \\
             122 & 231 & 231 \\
             123 & 321 & 321
        \end{tabular}
        \caption{Table of function values of $f_{1,3}$ and $f_{2,3}$ on words}
        \label{table_f}
    \end{table}
    
\begin{theorem}\label{thm:f_{1,n}f_{2,n}}
The fixed points of $f_{1,n}\circ f_{2,n}^{-1}$ are the involutions of $S_n$ but backwards.

\begin{proof}
Suppose $\sigma$ is a fixed point and consider a word $w = f^{-1}_{1,n}(\sigma) =f^{-1}_{2,n}(\sigma)$. Let us now suppose $f_{1,n}(w)[i] = f_{2,n}(w)[i] = n$ for some $i \in [n]$. Under the rules of $f_{1,n}$ this yields $w[n-i+1] = n-i+1$. Applying $f_{2,n}$, we must have $w[n] =n-i+1$. Furthermore  $w[j] < n-j+1$ for all $j \in [i+1, n]$ to ensure we don't place $n$ in the $f_{2,n}(w)$ before reaching position $j$.

It follows that $w[n-i+1] = w[n] = n-i+1$ if $f_{1,n}(w)[i] = f_{2,n}(w)[i] = n$. Let us now evaluate $f_{2,n}$ on this $w$. Since $w[n-i+1] = n-i+1$, we take $n-i+1$ and place it in position $n-i+1$ of $f_{2,n}(w)$ counting from the right among all numbers that are currently in the permutation. As there are only $n-i$ numbers we obtain $f_{2,n}(w)[1] = n-i+1$ at this step. Since $w[j] < n-j+1$ for all $j \in [i+1, n]$, in constructing $f_{2,n}(w)$, none of $[n-i+2, n]$ can be placed in front of $n-i+1$ which means $n-i+1$ will stay at position $1$.

Thus when we read this permutation backward, $n$ will be in position $n-i+1$ from the left while $n-i+1$ will be in position $n$. This is a $2$-cycle for $i > 1$ and a one cycle for $i = 1$. We can now remove $i$ and $n$ from his permutation and repeat the same reasoning on the remaining elements of the permutation to conclude that the fixed points backward is a self-inverse of $S_n$.

Conversely, let $\sigma$ be a permutation such that $\sigma$ backwards is an involution. Thus, we can suppose that when we write $\sigma$ backwards the resulting permutation maps $n$ to $k$ and $k$ back to $n$. Then $\sigma[1] = k, \sigma[n+1-k] = n$. Since $\sigma[n+1-k] = n$, we must have $f_{1,n}^{-1}(\sigma)[n] = k$. Similarly, since $\sigma[1]=k$, we must have $f_{2,n}^{-1}(\sigma)[k] = k$. Furthermore, $\forall j \in [k+1, n-1]$ we must have $f_{2,n}^{-1}(\sigma)[j] < j$ so that nothing gets inserted in front of $k$ when we apply $f_{2,n}$ to $f_{1,n}^{-1}(\sigma)$ to obtain $\sigma$.

We can now construct the permutation resulting from $f_{2,n}$ applied to the word $f_{1,n}^{-1}(\sigma)$. Since $f_{1,n}^{-1}(\sigma)[n] = k$, it follows that $f_{2,n}\circ f_{1,n}^{-1}(\sigma)[n] = k$. Furthermore, since only numbers less than $n$ will be added into the permutation $f_{2,n}\circ f_{1,n}^{-1}(\sigma)$ in positions $n$ through $n-k$ of $f_{1,n}^{-1}(\sigma)$. Thus, when we reach position $n-1+k$ of $f_{2,n}\circ f^{-1}(\sigma)$, the $k$-th smallest available number will be $n$ so $f_{2,n}\circ f_{1,n}^{-1}(\sigma)[n+1-k] = n$. Since there is no restriction on $k$ being different than $n$ in our reasoning above, it follows that the map $f_{2,n}\circ f_{1,n}^{-1}(\sigma)$ fixes the fixed point $n$ or transposition $(k, n)$ of $\sigma$ when written backwards. We can remove this pair and repeat the same reasoning on the remaining elements of $\sigma$ (treating them as a permutation on either $[n-2]$ or $[n-1]$ to obtain the desired result.
\end{proof}
\end{theorem}

Building on the previous theorem, we can generalize the $k$-cycles of $f_{1,n}\circ f_{2,n}^{-1}\circ f_{1,n}\circ f_{2,n}^{-1}$ with the following result. 

\begin{theorem}\label{thm:f12}
$f_{1,n}\circ f_{2,n}^{-1}\circ f_{1,n}\circ f_{2,n}^{-1}: S_n \to S_n$ is the identity. 
\begin{proof}
Suppose that $\sigma \in S_n$ was such that $\sigma[i] = n$. As with the previous theorem, We can apply the for maps in succession and deduce the corresponding position based on the rules given by $f_{1,n}$ and $f_{2,n}$. We get,
\begin{align}
    \sigma[i] = n &\implies f_{1,n}^{-1}(\sigma)[n] = n-i+1 \\
    &\implies f_{2,n}\circ f_{1,n}^{-1}(\sigma)[1] = n-i+1 \\
    &\implies f_{1,n}^{-1} \circ f_{2,n}\circ f_{1,n}^{-1}(\sigma)[n-i+1] = n-i+1 \text{ and } f_{1,n}^{-1} \circ f_{2,n}\circ f_{1,n}^{-1}(\sigma)[j]< j \text{ if } j > n-i+1 \\
    &\implies f_{2,n} \circ f_{1,n}^{-1} \circ f_{2,n}\circ f_{1,n}^{-1}(\sigma)[i] = n
\end{align}

We can now remove the last entry of $\sigma$, treat the remaining permutation as a permutation on $n-1$ and repeat the reasoning above to obtain the desired result. 
\end{proof}
\end{theorem}

\subsection{Fixed points of $f_{1,n}\circ f_{3,n}^{-1}$}

We characterize the fixed points of $f_{1,n}\circ f_{3,n}^{-1}$ with the next theorem. 
\begin{theorem}
The number of fixed points of the map $f_{1,n} \circ f_{3,n}^{-1}: S_n \to S_n$ is $2^{n-1}$. Namely, the fixed points are of the form $f_{1,n}(w) = f_{3,n}(w)$ where $w \in [1] \times [2] \times \dots \times [n]$ only has $1$ and $2$'s. 

\begin{proof}
Let us suppose that $w[n] = i$ and $w[n-1] = j$ and try to show that $i = 1, 2$. Suppose that $i> 1$, based on the definition of $f_{3,n}$, we have the first cycles of $f_{3,n}(w)$ as follows 
\begin{equation}
    f_{3,n}(w) = \begin{cases}
     (1i)(2\cdots & \text{if } i>2, j = 1\\
       (12)(3\cdots & \text{if } i=2, j = 1\\
    (1ij\cdots &\text{if } 1 <j < i\\
     (1i[j-1]\cdots &\text{if } j\geq i>1
    \end{cases}
\end{equation}

In other words, 

\begin{equation}
    f_{3,n}(w)[i] = \begin{cases}
    j &\text{if } j < i\\
    j-1 &\text{if } j\geq i
    \end{cases}
\end{equation}

On the other hand, evaluating $f_{1,n}$ to $w$ yields:
\begin{equation}
    f_{1,n}(w)[2] =\begin{cases}
    j &\text{if } j < i\\
    j-1 &\text{if } j\geq i
    \end{cases}
\end{equation}

We now use $f_{3,n}(w) = f_{1,n}(w)$ to conclude that $i = 2$ as desired. We remove $w[n]$ and repeat the same reasoning on the first $n -1$ letters of $w$, treating the resulting permutation as an element of $S_{n-1}$.

Conversely, suppose $w$ only has $1$'s and $2$'s. We can prove $f_{1,n}(w) = f_{3,n}(w)$ by induction on $n$ with the base case being trivial. Suppose $w[n] = 1$. Then $f_{1, n}(w)[1] = f_{3,n}(w)[1] = 1$. On the other hand, if $w[n] = 2$ then $f_{1,n}(w)[1] = f_{3,n}(w)[1] = 2$. We can now treat the remaining letters of the word as a subset of $[n-1]$, and apply the inductive hypothesis to complete the proof.
\end{proof}
\end{theorem}

\subsection{Fixed Points of $f_{2,n}\circ f_{4,n}^{-1}$}

We characterize the fixed points of $f_{1,n}\circ f_{4,n}^{-1}$ with the next theorem. 
\begin{theorem}
The number of fixed points of the map $f_{2,n} \circ f_{4,n}^{-1}: S_n \to S_n$ is $2$ for $n \geq 2$. Namely, the fixed points are of the form $f_{1,n}(w) = f_{4,n}(w)$ where $w$ are words such that $w[2] \in \{1, 2\},$ while $w[i] = 1$ for $i \neq 2$. 

\begin{proof}
Suppose $n > 2$, and let our corresponding word for the fixed point be $w$ such that $w[n] = i$ and $w[n-1] = j$. We hope to show that $w[n] = 1$. Suppose this wasn't the case, then we have the following 3 cases

\emph{Case 1:} $1< i < j$. In this case by the rules of $f_{2,n}$ we have $f_{2,n}(w)[n-j] = n-1$ and $f_{2,n}(w)[n-i+1] = n$. By the rules of $f_{4,n}$ we have $f_{4,n}(w)[j-1] = n-1$ while $f_{4,n}(w)[i-1] = n$. Since $f_{2,n}(w) = f_{4,n}(w)$, it follows that $n-j = j-1 \implies j = (n+1)/2$ while $n-i+1 = i-1 \implies i = (n+2)/2$. However, then one of $i$ or $j$ isn't an integer, which is impossible

\emph{Case 2:} $i = j>1$. In this case by the rules of $f_{2,n}$ we have $f_{2,n}(w)[n-j] = n-1$ and $f_{2,n}(w)[n-i+1] = n$. By the rules of $f_{4,n}$ we have $f_{4,n}(w)[n] = n-1$ while $f_{4,n}(w)[i-1] = n$. Since $f_{2,n}(w) = f_{4,n}(w)$, it follows that $j = 0$ which is impossible. 

\emph{Case 3:} $i > j > 1$, In this case by the rules of $f_{2,n}$ we have $f_{2,n}(w)[n-j+1] = n-1$ and $f_{2,n}(w)[n-i+1] = n$. By the rules of $f_{4,n}$ we have $f_{4,n}(w)[j-1] = n-1$ while $f_{4,n}(w)[i-1] = n$. Since $f_{2,n}(w) = f_{4,n}(w)$, it follows that $i = j = (n+2)/2$, which contradicts $i > j$. 

Conversely, suppose $w$ is a word such that $w[i] = 1$ for all $i \neq 2$ and $w[2] \in \{1, 2\}$. When $w = 11\dots 1$ we have $f_{2,n}(w) = 123\dots n$ and $f_{4,n}(w) = (1)(2)(3)\dots (n) = 123\dots n$. In the case $w = 121\dots 1$, we have $f_{2,n}(w) = 2134\dots n$ and $f_{4,n}(w) = (12)(3)(4)\dots (n) = 2134\dots n$. This concludes our proof. 
\end{proof}
\end{theorem}

\section*{Conclusions and Future Work}

\indent \indent In this paper we evaluated $4$ different types of bijective maps from words $[1]\times [2] \times \dots \times [n]$ to permutations $S_n$. We looked at the $6$ possible types of compositions arising from these maps that give bijections between $S_n$ and tried to characterize these maps. We were able to show $f_{1,n}\circ f_{2,n}^{-1}$ maps each permutation to its inverse written backwards. Furthermore, we completely characterized the fixed points of $f_{1,n}\circ f_{3,n}^{-1}$ and $f_{2,n}\circ f_{4,n}^{-1}$. For the remaining, we showed numerically in the appendix the complex cycle structure involved. For future work, it would be interesting to derive asymptotic bounds on these maps regarding the number of fixed points or on the order of the $k$-cycle structure. We also believe the following conjecture to be true.

\begin{conjecture}
There is only $1$ fixed point of $f_{4,n}\circ f_{1,n}^{-1}:S_n \to S_n$, namely it is the identity permutation and corresponds to the word $1...1$.
\end{conjecture}

\section*{Acknowledgements}
The author would like to thank Peter Kagey and Richard Arratia for their meaningful discussions and suggestions on this problem. 

	\newpage
	\bibliographystyle{plain}
	\bibliography{references}
	\newpage

\section{Appendix} 
\subsection{More on the Bijection $f_{1,n}\circ f_{2,n}^{-1}$}\label{appendix:f12}

\subsubsection{A Sibling Map to $f_{1, n}$ and $f_{2,n}$}
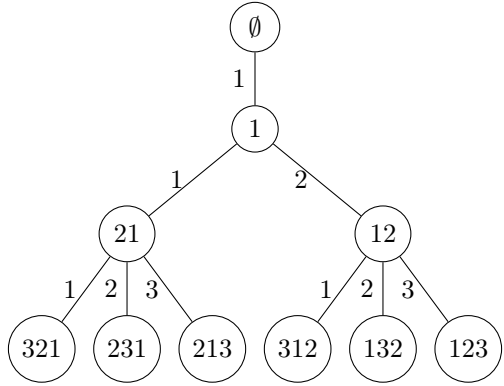
\begin{figure}
\centering
\begin{forest}
for tree={circle,draw, l sep=20pt}
[$\emptyset$ 
   [1, edge label={node[midway,left] {1}}
   [21, edge label={node[midway,left] {1}}
      [321, edge label={node[midway,left] {1}} ] 
      [231, edge label={node[midway,left] {2}}] 
      [213, edge label={node[midway,left] {3}}]
    ]
    [12, edge label={node[midway,left] {2}}
      [312, edge label={node[midway,left] {1}}] 
      [132, edge label={node[midway,left] {2}}] 
      [123, edge label={node[midway,left] {3}}]
    ]
  ]
]
\end{forest}
\caption{Tree diagram of $g_{1,n}$ for $n=3$.}
\label{tree_g_1}
\end{figure}

\begin{figure}
\centering
\begin{forest}
for tree={circle,draw, l sep=20pt}
[$\emptyset$ 
   [1, edge label={node[midway,left] {1}}
   [21, edge label={node[midway,left] {1}}
      [321, edge label={node[midway,left] {1}} ] 
    ]
    [31, edge label={node[midway,left] {2}}
      [231, edge label={node[midway,left] {1}}] 
    ]
  ]
  [2, edge label={node[midway,left] {2}}
   [12, edge label={node[midway,left] {1}}
      [312, edge label={node[midway,left] {1}} ] 
    ]
    [32, edge label={node[midway,left] {2}}
      [132, edge label={node[midway,left] {1}}] 
    ]
  ]
  [3, edge label={node[midway,left] {3}}
   [13, edge label={node[midway,left] {1}}
      [213, edge label={node[midway,left] {1}} ] 
    ]
    [23, edge label={node[midway,left] {2}}
      [123, edge label={node[midway,left] {1}}] 
    ]
  ]
]
\end{forest}
\caption{Tree diagram of $g_{2,n}$ for $n=3$.}
\label{tree_g_2}
\end{figure}
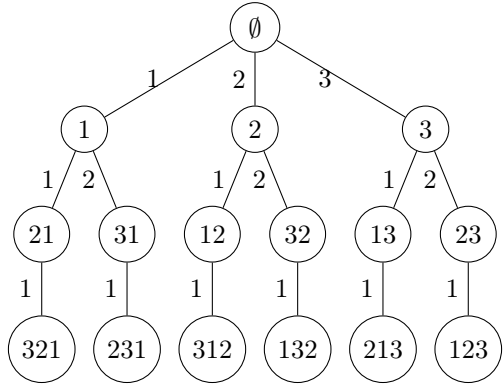

We showed that the map $f_{1,n}\circ f_{2, n}^{-1}:S_n \to S_n$ maps each permutation to its inverse but written backwards. In hopes of mapping each permutation to its true inverse, let us now 'flip' the maps $f_{1,n}$ and $f_{2,n}$ to obtain new maps $g_{1,n}, g_{2,n}: [1]\times \cdots \times [n] \to S_n$. More explicitly, the maps $g_{1,n}, g_{2,n}$ are constructed the same as $f_{1,n}$ and $f_{2,n}$, at each step instead of inserting from the right, we insert from the left. We draw out the tree diagrams for $g_{1,n}$ and $g_{2,n}$ for $n = 3$ in Figures \ref{tree_g_1} and \ref{tree_g_2} respectively.

\subsubsection{Fixed Points of $g_{1,n}\circ g_{2,n}^{-1}$}

 We take the definitions of $g_{1,n}$ and $g_{2,n}$ from the previous section, and consider $g_1\circ g_2^{-1}: S_n \to S_n$. If we tabulate the results from the tree diagrams in Figure \ref{tree_g_1} and Figure \ref{tree_g_2} as shown in Table \ref{table_g}, then looking at the middle column of this table to the rightmost column, we can conjecture that $g_1\circ g_2^{-1}$ maps a permutation to its inverse. We prove this with the next theorem.

\begin{table}[]
        \centering
        \begin{tabular}{c|c|c}
            word & $g_{2,n}$ & $g_{1,n}$ \\
            \hline
             111 & 321 & 321 \\
             112 & 312 & 231 \\
             113 & 213 & 213 \\
             121 & 231 & 312 \\
             122 & 132 & 132 \\
             123 & 123 & 123
        \end{tabular}
        \caption{Table of function value of $g_{1,n}$ and $g_{2,n}$ on words}
        \label{table_g}
    \end{table}

\begin{theorem}\label{thm:fixed}
The fixed points of $g_{1,n}\circ g_{2,n}^{-1}$ are the involutions of $S_n$
\begin{proof}

Suppose $\sigma$ is a fixed point and consider a word $w = g^{-1}_{1,n}(\sigma) = g^{-1}_{1,n}(\sigma)$. Since $f_{1,n}(w)$ that $g_{1,n}(w)[i] = g_{2,n}(w)[i] = n$. Under the rules of $g_{1,n}$ the number at position $k$ of $w$ determines where $k$ goes in $g_{1,n}(w)$. Thus, in order for $g_{1,n}(w)[i] = n$ to hold we must have $w[n]=i$. Under the rules of the map $g_{2,n}$, the number at position $w[k]$ determines what $g_{2,n}(w)[k]$ is. Thus, in order for $g_{2,n}(w)[i] = n$ to hold we must have $w[i] = i$. Furthermore, $w[j] < j$ for all $j \in [n-i+1, n]$ to ensure we don't place $n$ in the $g_{2,n}(w)$ before reaching position $j$. 

Using $w[i]=w[n] = i$, let us now evaluate $f_{1,n}$ on this $w$. Since $w[i] = n-i+1$ means we take $i$ and place it in position $n-i+1$ of $g_{1,n}(w)$ counting from the left among all numbers that are currently in the permutation. Since there are only $i$ numbers this means that $g_{1,n}(w)[i] = i$. We note that by the remark above, in constructing $g_{1,n}(w)$, none of $[n-i+1, n]$ can be placed behind $i$ which means $i$ will stay at the back. By construction, $g_{1,n}(w)[i] = n$. Thus, this permutation will have the transposition $(i, n)$ and in the case that $i = n$, $n$ will a fixed point. We can now remove $i$ and $n$ from his permutation and repeat the same reasoning on the remaining elements of the permutation to conclude that the fixed points are involutions of $S_n$.

Conversely, let $\sigma$ be a permutation such that $\sigma$ is an involution. Then we can suppose that $\sigma$ maps $n$ to $k$ and $k$ back to $n$ so that $\sigma[k] = n, \sigma[n] = k$. They key here is to consider $g_{2,n}\circ g_{1,n}^{-1}$ instead of $g_{1,n}\circ g_{2,n}^{-1}$. As before, $g_{1,n}^{-1}(\sigma)[n]$ determines the position of $n$ in $\sigma$ counting from the left, so since $\sigma[k] = n$, we must have $g_{1,n}^{-1}(\sigma)[n] = k$. Similarly, since $\sigma[n]=k$, we must have $g_{1,n}^{-1}(\sigma)[k] = k$. Furthermore, $\forall j \in [k+1, n-1]$ we must have $g_{1,n}^{-1}(\sigma)[j] < j$ so that nothing gets inserted in behind $k$ when we apply $g_{1,n}$ to $g_{1,n}^{-1}(\sigma)$.

We can now construct the permutation from $g_{2,n}$ applied to the word $w = g_{1,n}^{-1}(\sigma)$. Furthermore, only numbers less than $n$ will be added into the permutation $g_{2,n}\circ g_{1,n}^{-1}(\sigma)$ before reaching position $k$ in the word $g_{1,n}^{-1}(\sigma)$. Thus, when we reach position $k$ of $g^{-1}(\sigma)$, the $k$-th smallest available number will be $n$ so $g_{2,n}\circ g_{1,n}^{-1}(\sigma)[k] = n$. Since there is no restriction on $k$ being different than $n$ in our reasoning above, it follows that the map $f_{2,n}\circ f_{1,n}^{-1}(\sigma)$ fixes the fixed point $n$ or transposition $(k, n)$ of $\sigma$. We can remove this pair and repeat the same reasoning on the remaining elements of $\sigma$ to obtain the desired result.
\end{proof}
\end{theorem}

\subsubsection{A New Sequence}
Let us consider the map denoted $h_{1,n}, h_{2,n}: [n!] \to S_n$ which follows the same construction as $g_{1,n}, g_{2,n}$ but instead of it mapping $S_n \to S_n$ we consider this map $[n] \to [n]$. We take the exact same construction of $g_{1,n}, g_{2,n}$ but associate each word in lexicographic order to an integer. This is equivalent numbering each leaf of the trees given by Figure \eqref{tree_g_1} and Figure \eqref{tree_g_2} and numbering that from left to right. We list out the explicit values of $h_{1,n}: [3!] \to S_3$ as given in Table \ref{table_h_2}. Similarly, for $h_{2,n}$, we have Table \ref{table_h_2}.

\begin{table}[]
        \centering
        \begin{tabular}{c|c}
            $k$ & $h_{3,n}(k)$ \\
            \hline
             1 & 321 \\
             2 & 312 \\
             3 & 213 \\
             4 & 231 \\
            5 & 132\\
             6 & 123
        \end{tabular}
        \caption{Table of function value of $h_{1,n}: [3!] \to S_3$ for $n = 3$}
         \label{table_h_1}
    \end{table}

\begin{table}[]
        \centering
        \begin{tabular}{c|c}
            $k$ & $g_3(k)$ \\
            \hline
             1 & 321 \\
             2 & 231 \\
             3 & 213 \\
            4 & 312 \\
            5 & 132 \\
            6 & 123
        \end{tabular}
        \caption{Table of function value of $h_{2,3}: [3!] \to S_3$}
        \label{table_h_2}
    \end{table}
    
We now consider the fixed points of $h_{1,n}\circ h_{2,n}^{-1}$ as we have done before and the results are starting from $n=2$ are as follows;  $2, 3, 3, 3, 10, 5, 4, 5,  13, 3, 6, 5$. It is currently sequence A347208 \footnote{\url{https://oeis.org/A347208}} published by the author.

\subsection{Numerical Results for Cycle Structures of Bijectiions $S_n \to S_n$}\label{appendx:numerical}
In Theorem \ref{thm:f12} we completely characterized the map $f_{1,n} \circ f_{2,n}^{-1}$. For the other maps, $f_{1,n}\circ f_{2,n}^{-1}$, $f_{1,n}\circ f_{3,n}^{-1}$, $f_{1,n}\circ f_{4,n}^{-1}$, $f_{2,n}\circ f_{3,n}^{-1}$, $f_{2,n}\circ f_{4,n}^{-1}$, and $f_{3,n}\circ f_{4,n}^{-1}$, we have numerically evaluated the number of $k$ cycles for small values of $n$ in Tables \ref{table_f13}, \ref{table_f14}, \ref{table_f23}, \ref{table_f24}, and \ref{table_f34} respectively. They are written in array notation with the $k$-th element being the number of $k$ cycles of this map. Where there is an entry of the form 0: m for some positive integer m, this means there are m consecutive 0's in the array. Since these maps have permutations with large k-cycles, such notation is used to make the k cycle arrays more compact.
\begin{table}[]
        \centering
        \begin{tabular}{@{}p{1cm}@{}p{15cm}}
            $n$ & Cycle Structure \\
            \hline
             1 & [1] \\
             2 & [2] \\
             3 & [4, 2] \\
             4 & [8, 6, 0, 4, 0, 6] \\
             5 & [16, 16, 12, 28, 0, 18, 0, 8, 0, 10, 0, 12] \\
             6 & [32, 44, 36, 84, 0, 48, 0, 24, 18, 30, 0, 36, 0, 14, 0: 28, 86, 88, 0: 3, 96, 0: 35, 84] \\
             7 & [64, 120, 102, 244, 0, 156, 14, 64, 54, 100, 0, 96, 26, 56, 0, 16, 0: 5, 22, 0: 3, 52, 54, 56, 0: 7, 72, 0: 6, 258, 264, 0: 3, 288, 0: 11, 60, 0: 5, 132, 0: 9, 76, 0: 3, 160, 0: 3, 420, 0: 6, 182, 92, 0: 75, 168, 0: 57, 226, 0: 11, 476, 0: 5, 488, 0: 137, 382]
        \end{tabular}
        \caption{Table of the cycle structures of the map $f_{1,n}\circ f_{3,n}^{-1}: S_n \to S_n$.}
        \label{table_f13}
    \end{table}
    \begin{table}[]
        \centering
        \begin{tabular}{@{}p{1cm}@{}p{15cm}}
            $n$ & Cycle Structure \\
            \hline
             1 & [1] \\
             2 & [2] \\
             3 & [1, 0: 3, 5] \\
             4 & [1, 2, 0, 4, 0: 12, 17] \\
             5 & [1, 2, 0: 7, 20, 0: 7, 18, 0: 60, 79]\\
             6 & [1, 0: 2, 4, 5, 0: 704, 710] \\
             7 & [1, 2, 0: 13, 16, 0: 191, 208, 0: 276, 485, 0: 66, 552, 0: 1297, 1850, 0: 75, 1926]
        \end{tabular}
        \caption{Table of the cycle structures of the map $f_{1,n}\circ f_{4,n}^{-1}: S_n \to S_n$.}
        \label{table_f14}
    \end{table}

        \begin{table}[]
        \centering
        \begin{tabular}{@{}p{1cm}@{}p{15cm}}
            $n$ & Cycle Structure \\
            \hline
             1 & [1] \\
             2 & [2] \\
             3 & [2, 4] \\
             4 & [5, 8, 6, 0, 5] \\
             5 & [4, 24, 3, 0: 2, 18, 0: 2, 9, 10, 0: 10, 21, 0: 9, 31]\\
             6 & [12, 54, 9, 0, 5, 54, 0: 3, 20, 0: 3, 14, 0: 3, 18, 0: 20, 39, 0: 2, 42, 0: 4, 47, 0: 10, 58, 0: 3, 62, 0: 223, 286] \\
             7 & [11, 140, 3, 0, 5, 174, 7, 0, 9, 80, 0: 2, 13, 28, 0: 3, 36, 0: 2, 21, 0: 9, 31, 0: 10, 84, 0: 15, 116, 0: 3, 124, 0: 15, 78, 0: 15, 94, 0: 148, 243, 0: 42, 572, 0: 183, 470, 0: 2230, 2701]
        \end{tabular}
        \caption{Table of the cycle structures of the map $f_{2,n}\circ f_{3,n}^{-1}: S_n \to S_n$.}
        \label{table_f23}
    \end{table}
            \begin{table}[]
        \centering
        \begin{tabular}{@{}p{1cm}@{}p{15cm}}
            $n$ & Cycle Structure \\
            \hline
             1 & [1] \\
             2 & [2] \\
             3 & [2, 0: 2, 4] \\
             4 & [2, 0, 6, 4, 0, 12] \\
             5 &[2, 0, 6, 4, 0, 12, 0: 17, 96]\\
             6 & [2, 0, 6, 4, 20, 12, 0: 3, 40, 0: 9, 60, 0: 3, 96, 0: 5, 180, 0: 29, 120, 0: 29, 180] \\
             7 & [2, 0, 6, 4, 20, 12, 0: 3, 40, 0, 48, 0: 7, 60, 0: 3, 192, 0: 5, 180, 0: 5, 216, 0: 11, 624, 0: 11, 120, 0: 11, 720, 0: 17, 180, 0: 5, 672, 0: 11, 648, 0: 35, 864, 0: 71, 432]
        \end{tabular}
        \caption{Table of the cycle structures of the map $f_{2,n}\circ f_{4,n}^{-1}: S_n \to S_n$.}
        \label{table_f24}
    \end{table}
                \begin{table}[]
        \centering
        \begin{tabular}{@{}p{1cm}@{}p{15cm}}
            $n$ & Cycle Structure \\
            \hline
             1 & [1] \\
             2 & [2] \\
             3 & [1, 2, 3] \\
             4 & [2, 2, 9, 0: 7, 11] \\
             5 &[3, 2, 9, 0: 2, 12, 14, 0: 2, 10, 0, 12, 0: 10, 23, 0: 11, 35]\\
             6 & [4, 4, 27, 0, 10, 24, 7, 0, 9, 0: 7, 17, 0: 9, 27, 0: 13, 41, 0: 4, 46, 0: 28, 75, 0, 77, 0: 7, 85, 0: 4, 90, 0: 86, 177]\\
             7 & [7, 6, 36, 0, 10, 96, 0: 2, 9, 0: 9, 19, 40, 0: 3, 24, 25, 26, 0: 15, 42, 0: 18, 61, 0: 23, 85, 0: 24, 110, 0, 112, 0: 10, 123, 124, 0: 87, 212, 0: 3, 216, 0: 43, 260, 0: 16, 277, 0: 27, 305, 0: 85, 391, 0: 415, 807, 0: 809, 1617]
        \end{tabular}
        \caption{Table of the cycle structures of the map $f_{3,n}\circ f_{4,n}^{-1}: S_n \to S_n$.}
        \label{table_f34}
    \end{table}

\end{document}